\documentclass[12pt]{article}
\usepackage{epsfig}

\makeatletter
\renewcommand{\@biblabel}[1]
\makeatother \textwidth160mm \textheight250mm
\newtheorem{Sa}{\sc   Proposition}
\newtheorem{Le}[Sa]{\sc  Lemma}

\newcommand{\RR}{{\rm I\kern-0.14em R }}
\newcommand{\NN}{{\rm l\kern-0.14em N}}
\newcommand{\CC}{{\rm\raise 0.192ex\vbox{\hrule height
                  1.22ex width 0.8pt}\kern-0.29em C}}
\newcommand{\ZZ}{ {\sf Z}\hspace{-0.4em}{\sf Z}\ }

\begin{document}

\begin{center} \bf \Large
The group law for the Jacobi variety \\
of plane curves
\end{center}

\vspace{0.1cm}

\begin{center}
{\large Frank Leitenberger}
\end{center}

\begin{center} \sl \small
Fachbereich Mathematik, Universit\"at Rostock, Rostock, D-18051, Germany, \\
 frank.leitenberger@uni-rostock.de
\end{center}

\vspace{0.1cm}

%\begin{flushleft}
%\small \bf Abstract
%\end{flushleft}

\begin{flushleft} \small
We extend the group law of curves of degree three by chords and
tangents to the Jacobi variety of plane curves of degree $n\geq 4$
by replacing points by point groups and lines by algebraic curves.
The curves are nonsingular or have simple singularities. In the
cases $n=4,5,6$ we have a close analogy to the case $n=3$. We describe
an algorithm using Groebner bases.
\end{flushleft}

\begin{flushleft}
\end{flushleft}

\section{Introduction}

\noindent
The Jacobi variety of an algebraic curve of genus $g$
plays a central role in algebraic geometry. Explicit descriptions
of the group law play a less important role in the history of the
subject.
%Mazur remarked 1986: "{\it ... a naive attempt to
%generalize this group structure {\rm [of degree 3 curves]} to
%curves of higher degree (even quartics) will not work.}", cf.
%\cite{Ma}, p. 230.
With the development of cryptography
algorithms arose for the group law. In 1987 D.G. Cantor described the
group law of a hyperelliptic curve in the context of cryptography
in analogy to the Gau{\ss} composition of quadratic forms,
cf. \cite{Ca,K}.

Later group laws of other classes of plane curves
were described.
The papers \cite{Ga, Har, He} are based on the analogy
of the Jacobian group with class groups in number theory.
Their methods are restricted to curves
with a special type of infinite points.
Other papers are concerned with special types of curves
(e.g. Picard curves \cite{Es,FO}, $C_{ab}$-curves \cite{Ar}).
The papers \cite{Vo,Hu} describe general ideas for an algorithm
for arbitrary plane curves. They are based on the theory of
plane curves. These algorithms are not practical from a computational
point of view.

The content of this paper is an algorithm
for smooth plane curves or of plane curves with simple double points.
We give an elementary presentation using projective curves,
Riemann surfaces and commutative algebra.
We consider arbitrary
complex nonsingular plane $n$-curves $C$ $(n\geq 4)$ with an
arbitrary zero point. We describe a geometric group law
intersecting $C$ with certain algebraic $m$-curves. For $n=4,5,6$
we intersect $C$ with $(n-1)$-curves and we have a very close
analogy to the group law of an elliptic curve.
On the basis of the
geometric group law we give an algebraic description. We represent
divisors by certain homogeneous ideals and describe the group law
by ideal operations.
All operations are rational and can be
carried out by Groebner basis operations.
Infinite points do not play a special role.
The construction works also for curves with simple double points.
We give an example for a curve of genus 2.
In the last sections we discuss hyperelliptic curves, Picard curves
and the case of fields of characteristic $p$.

Because of the use of Groebner bases it is difficult to give a
realistic complexity analysis.
We used the computer algebra program Mathematica 4.0
forming the computations.

\section{Nonsingular curves and intersections}

\noindent Consider a nonsingular projective plane curve $C$ of
degree $n\geq 4$ defined by the irreducible polynomial
\[
F(x,y,z):= \sum_{i+j\leq n} a_{ij}\ x^iy^jz^{n-i-j}.
\]
One can consider $C$ as a Riemann surface of genus
$g=\frac{(n-1)(n-2)}{2}$ (i.e. $g=3,6,10,\ \cdots $). In every
point one of the coordinates $x,y,z$ defines a complex chart. We
have a $g$-dimensional space of holomorhic differentials.
Every holomorphic differential admits an explicit
representation in the form
\[
G(x,y,z) \frac{ \left| \begin{array}{ccc}
\gamma_1 & \gamma_2 & \gamma_3 \\
x & y & z \\
dx & dy & dz
\end{array} \right|}{\gamma_1F_x+\gamma_2F_y+\gamma_3F_z}
\]
where $G(x,y,z)$ is a homogeneous polynomial of degree $n-3$ and
$\gamma_1,\gamma_2,\gamma_3$ are three complex numbers chosen depending on $x,y,z$
such that the denominator is not zero
(cf. \cite{Cl}). Furthermore let
$H(P_1+\ \cdots\ +P_t)$
be the subspace of holomorphic differentials whose zero divisors
contain the points $P_1,\ \cdots\ ,P_t$.

Now we consider $m$-curves $C_m$. It is allowed that $C_m$ has
multiple irreducible components. First we consider the case $1\leq
m \leq n-1$. $C_m$ has $\frac{(m+1)(m+2)}{2}$ parameters. Given
$b_m:=\frac{m(m+3)}{2}$ points on $C$ we have a (not necessary
unique) $C_m$-curve through these points. (A multiple intersection
point will correspond to a multiple contact of $C_m$ with $C$).
Because of the irreducibility of $C$ an $m$-curve $C_m$ has
$mn$ intersections with $C$. Therefore we have
\[
c_m:=mn-\frac{m(m+3)}{2}=\frac{m(2n-m-3)}{2} =g-\left(
\begin{array}{c} n-m-1 \\ 2
\end{array} \right)
\]
further intersections.

For $n=4,5,6$ we will be interested in the cases $m=n-1,n-2$ and
$n-3$. Here we have $c_{n-1}=g, c_{n-2}=g, c_{n-3}=g-1$ and
$b_{n-1}=g+2n-2, b_{n-2}=g+n-2, b_{n-3}=g-1$.

For $n>6$ we will need also curves $C_m$ of degree $m\geq n$. In
this case we have the difficulty that it is possible, that $C$ is
an irreducible component of $C_m$. For $m\geq n$ we fix a monomial
$x^ky^lz^{n-k-l}$ of $F$ with $a_{kl}\neq 0$. Then we consider the
$m$-curves
\[
\sum_{
\begin{array}{c} (i,j)\ \ {\rm with }
  \\ x^ky^lz^{n-k-l}\not{\kern-0.14em\mid
} \ x^iy^jz^{m-i-j}
\end{array} }
\beta_{ i,j,m-i-j}\ x^iy^jz^{m-i-j}\ .
\]
They form a linear system of dimension {\footnotesize $\left(
\begin{array}{c} m+2 \\ 2
\end{array} \right)-
\left( \begin{array}{c} m-n+2 \\ 2
\end{array} \right)$} $-1=mn-g$
of curves without common components with $C$. Given $b_m:=mn-g$
points on $C$ we find a $C_m$-curve of this special form through
these points. The curves $C$ and $C_m$ have no common components
and therefore they have $mn$ intersections, i.e. we have $c_m=g$
further intersections.

\section{Jacobi varieties and Reduced Divisors}

The Jacobi variety of $C$ is the Abelian group
\[ Jac(C)=Div^0(C)/Div^P(C).\]
Here $Div^0(C)$ denotes the group of divisors of degree $0$ and
$Div^P(C)$ is the subgroup of principal divisors (i.e. the zeros
and poles of analytic functions), cf. \cite{Ha}. If $G(x,y,z)$
defines an $m$-curve $C_m$, then $G(x,y,1)$ defines a meromorphic
function with zeros at the finite intersections of $C_m$ and with
poles of order $m-h$ at the infinite points $P$ of $C$ if $P$ is a
$h$-fold intersection with $C_m$. Let $D_\infty $ be the divisor
of the infinite points of $C$. Then we have $P_1+\ \cdots\ +P_{mn}
\sim mD_\infty  $, if $P_1,\ \cdots\ ,P_{mn}$ are the (possibly
infinite) intersections of $C_m$ and $C$.

Wie fix an arbitrary point $P_0$ on $C$.
We call a divisor $D$ of the form
$D=P_1+\ \cdots\  +P_t-tP_0$ with $P_1,\ \cdots\ ,P_t\neq P_0$
a semireduced divisor.
We call this semireduced divisor $D$ a reduced divisor
if there is no divisor
$D'={P'}_1+\ \cdots\  +{P'}_s-sP_0$
with ${P'}_1,\ \cdots\ ,{P'}_s\neq P_0$,
$D'\sim D$ and $s< t$.

\begin{Sa}
We find in every divisor class of
Jac(C) a unique reduced divisor with $t\leq g$.
\end{Sa}

\noindent
{\it Proof:}
{\it Existence:} Let $D\in Div^0(C)$.
By the Riemann Roch theorem we have \\
$\dim L(D+gP_0)\geq g-g+1=1$. It follows
\[ D+gP_0+(f)=P_1+\cdots +P_g \]
with $(f)\in L(D+gP_0)$ and certain $P_i\in C$.
Therefore we have $D\sim P_1+\cdots +P_g -gP_0$, i.e
the set of semireduced divisors with $t\leq g$
is not empty and we can find reduced divisors.

\noindent
{\it Uniqueness:}
Let $D =P_1+\cdots +P_t -tP_0\sim
     D'=Q_1+\cdots +Q_t -tP_0$
be two different reduced divisors with $P_i,Q_i\neq P_0$.
There is a function $f$ with $(f)=D-D'$.
Then the divisor
$D'+(f-f(P_0))$ has the form
$R_1+\cdots R_s -sR_0$ with $s<t$.
This is a contradiction.
$\bullet$

\noindent
{\it Remark:} One can choose $P_0$ arbitrarily,
but the structure of a reduced divisor might vary with this choice.
For the reason of simplicity it can be useful
to choose for $P_0$ an exceptional point
(flex, rational point, infinite point etc.).

\noindent
{\it Remark:}
We mention that not all semireduced divisors with $t\leq g$
are reduced. In the generic case
we have $t=g$ for a reduced divisor.

%A large class of reduced divisors is given by
%\begin{Sa}
%The divisor $D=P_1+\ \cdots\ +P_g-gP_0$ with $P_1,\ \cdots\
%,P_g\neq P_0$ is reduced if $P_1,\ \cdots\ ,P_g$ lie not on a
%projective $(n-3)$-curve.
%\end{Sa}

%\noindent {\it Proof.} Let $D$ be reduced and $P_1,\ \cdots\ ,P_g$
%lie on a projective $(n-3)$-curve $C_{n-3}$. Then $C_{n-3}$ has
%with $C$ $n(n-3)-\frac{(n-1)(n-2)}{2}=g-2$ further intersection
%points. They determine with $P_0$ a further (not necessary unique)
%$(n-3)$-curve ${C'}_{n-3}$. ${C'}_{n-3}$ has with $C$ $g-1$
%further intersections $Q_1 , \ \cdots\ ,Q_{g-1}$. It follows
%$P_1+\ \cdots\ +P_g \sim  Q_1+\ \cdots\ +Q_{g-1}+P_0$. Therefore
%$D$ is not reduced. $\bullet$

We have the following abstract characterisation of reduced
divisors.
\begin{Sa}
The following three assertions are equivalent.

\noindent (i) The divisor $D=P_1+\ \cdots\ +P_t-tP_0$ with $t\leq
g$ and $P_1,\ \cdots\ ,P_t\neq P_0$ is reduced.

\noindent (ii) $\dim\ H(P_1+\ \cdots\ +P_t) = g-t$.

\noindent (iii) The dimension of the linear system of
$(n-3)$-curves, which vanish on $P_1,\ \cdots\ ,P_t$
does not exceed $g-t$.
\end{Sa}

\noindent {\it Proof.} The proposition is a consequence of the
Riemann Roch theorem. $\bullet$

For $n=4,5,6$ we have an explicit description of the set of
reduced divisors.

\begin{Sa} (i) Let $n=4$.
The divisor $D=P_1+\ \cdots\ +P_t-tP_0$ with $t\leq g=3$ and
$P_1,\ \cdots\ ,P_t\neq P_0$ is reduced if and only if $t=3$ and
$P_1,P_2 ,P_3$ lie not on a projective line or $t<3$.

\noindent (ii) Let $n=5$. The divisor $D=P_1+\ \cdots\ +P_t-tP_0$
with $t\leq g=6$ and $P_1,\ \cdots\ ,P_t\neq P_0$ is reduced if
and only if no 6 points lie on a conic section and no 4 points lie
on a line.

\noindent (iii) Let $n=6$. The divisor $D=P_1+\ \cdots\ +P_t-tP_0$
with $t\leq g=10$ and $P_1,\ \cdots\ ,P_t\neq P_0$ is reduced if
and only if no 10 points lie on a cubic and no 9 points form the
intersection of two cubics and no 8 points lie on a conic section
and no 5 points lie on a line.
\end{Sa}

\noindent {\it Proof.} $n\leq 2d+2$ points of the projective plane
fail to impose independent conditions on curves of degree $d$ if
and only if either $d+2$ of the points are collinear or $n=2d+2$
and the $n$ points lie on a conic, cf. \cite{Ei}. The proposition
follows from this fact, from the Cayley-Bacharach theorem, cf.
also \cite{Ei} and from Proposition 2.$\bullet$

\noindent
{\it Remark:}
An explicit  characterisation of reduced divisors for arbitrary
$n$ is associated to the classification problem of special
divisors on Riemann surfaces or to the Cayley-Bacharach
conjectures for algebraic curves. However there is an algorithm to
determine the reduced divisor (cf. below).

\section{The construction of the reduced divisor}

Let
\[
D=D^+ - D^- =P_1+\ \cdots\ P_s-Q_1-\ \cdots\  -Q_s
\]
with $P_i,Q_i\in C$ be an arbitrary divisor of degree zero.

At first we consider an $m$-curve with a polynomial $G(x,y,z)\neq
0\ {\rm mod}\ F$ with a minimal $m\geq n-2$ such that
\[
s+g\leq b_m:=\left\{  \begin{array}{ccc}
\frac{m(m+3)}{2} & {\rm for} & m<n \\
mn-g & {\rm for} & m\geq n
\end{array} \right.\
\]
through the $s$ points of $D^+$ and $(b_m-s)P_0$.
Because of our discussion at the end of section 2
the polynomial $G(x,y,z)\neq 0\ {\rm mod}\ F$ exists and
there are $g$ remaining intersections $R_1,\ \cdots\ ,R_g$.
We have
\[
D^+ +(b_m -s )P_0+R_1+\ \cdots\ +R_g \sim m D_\infty .
\]
Then we consider another $m$-curve $G'(x,y,z)\neq 0 \rm\ mod\ \it
F$ through the $s+g$ points of $D^-+R_1+\ \cdots\ +R_g$ and
$(b_m-s-g)P_0$. Because this curve is not necessary unique we
require a maximal additional contact $\alpha$ at $P_0$. Let $S_1,\
\cdots\ ,S_{g-\alpha}$ be the remaining intersections not equal to
$P_0$. We have
\[
D^-+R_1+\ \cdots\ +R_g +(b_m -s-g+\alpha )P_0 +S_1+\ \cdots\
+S_{g-\alpha} \sim m D_\infty .
\]
It follows
\[
D^+ - D^-  \sim S_1+\ \cdots\ +S_{g-\alpha} - (g-\alpha)P_0 .
\]

\begin{Sa}
\[
\overline{D}:=S_1+\ \cdots\ +S_{g-\alpha} - (g-\alpha)P_0
\]
is the reduced divisor for $D$.
\end{Sa}

\noindent {\it Proof.} From the above relations it follows
$D\sim\overline{D}$. Furthermore we consider the divisor
\[
D_1:=D^-+R_1+\ \cdots\ +R_g +(b_m -s-g+\alpha )P_0
-(m-n+3)D_\infty .
\]

\begin{Le}
All differentials of $H(D_1)$ have the form
\[
\frac{G_m(x,y,z)}{z^{m-n+3}}\ \frac{ \left| \begin{array}{ccc}
\gamma_1 & \gamma_2 & \gamma_3 \\
x & y & z \\
dx & dy & dz
\end{array} \right|}{(\gamma_1F_x+\gamma_2F_y+\gamma_3F_z)}
\]
with ${\rm deg}\ G_m=m$.
\end{Le}

\noindent {\it Proof.} Case 1 (the positive part of $D_1$ does not
contain infinite points): Meromorphic differentials are unique
determined by its principal part up to a holomorphic differential.
The space of principal parts is $(m-n+3)n$-dimensional. Therefore
the space of differentials with poles at most at $(m-n+3)D_\infty$
has the dimension $(m-n+3)n+g=mn-g+1$. Otherwise the space of
degree $m$ polynomials $G\ {\rm mod}\ F$ has the dimension
$mn-g+1$ for $m\geq n-3$ (cf. above). Therefore every differential
with poles at most at $(m-n+3)D_\infty$ has the above form.

Case 2 ($\alpha $ points of the positive part $D_1^+$ and the
negative part $D_1^-$ of $D_1$ cancel): Analogously the space of
differentials with negative part at most at $D_1^-$ has the
dimension $mn-g+1-\alpha$. On the other side the space of degree
$m$ polynomials $G {\rm mod}\ F$ with zeros at the $\alpha $
common points of $D_1^+$ and $D_1^-$ has also the dimension
$mn-g+1-\alpha$.$\bullet$

For $G_m=G'$ we obtain a differential in $H(D_1)$. Because $C'$
has maximal contact at $P_0$ with $C$ we have no further
differential in $H(D_1)$. Therefore we have
\[
\dim\ H(D_1 )=1.
\]

Because $D_1+S_1+\ \cdots\ +S_{g-\alpha }$ is the divisor of poles
and zeros of the differential with $G_m=G'$ it is a canonical
divisor. By the Riemann Roch theorem it follows
\[
\dim\ H(S_1+\ \cdots\ +S_{g-\alpha })        = \dim\ L(D_1) =
\]
\[
(s+g+b_m-s-g+\alpha-(m-n+3)n)-g+1+1\ = \]
\[  mn -(m-n+3)n -2g+\alpha
+2=\alpha.
\]
Therefore $\overline{D}$ is reduced according to Proposition 2.
$\bullet$

\section{The group law }

\noindent We represent the elements of the Jacobian by reduced
divisors. For a description of the group law it is sufficient to
reduce $-D$ and $D_1+D_2$ if $D,D_1,D_2$ are reduced divisors.

\noindent 1. {\it The inverse divisor}. Let $D=P_1+\ \cdots\
+P_g-gP_0$ where $P_1,\ \cdots\ ,P_g = P_0$ is allowed. In this
case we apply the above construction to $D^+=gP_0$ and $D^-=P_1+\
\cdots\ +P_g$.

\noindent 2. {\it The addition of reduced divisors}. Let
$D_1=P_1+\ \cdots\ P_g-gP_0$ and $D_2=P_{g+1}+\ \cdots\
+P_{2g}-gP_0$. In this case we apply the above construction to
$D^+=P_1+\ \cdots\ +P_{2g}$ and $D^-=2gP_0$.

For $n=4,5,6$ we can carry out the construction with
$(n-1)$-curves. Therefore we have a close analogy to the group law
of cubic curves. In these cases we have $b_{n-1} = g+2n-2 \geq
2g$. Instead of the $(n-2)$-curves of the cubic case we consider
$(n-1)$-curves. We construct an $(n-1)$-curve through $P_1,\
\cdots\ ,P_{2g},(2n-2-g)P_0$. We have the remaining points $R_1,\
\cdots\ ,R_{g}$. Then we construct an $(n-1)$-curve through $R_1,\
\cdots\ ,R_{g}$ with highest contact at $P_0$ and obtain the
remaining points $S_1,\ \cdots\ ,S_{t}$. Then $S_1+\ \cdots\
+S_{t}-tP_0$ is the reduced divisor for $D_1+D_2$. Analogously,
for $n\geq 7$ we carry out the construction with $m$-curves
$G(x,y,z)\neq 0\ {\rm mod}\ F$ with $m\geq n$ and $b_m\geq 2g$
(i.e. $mn\geq 3g$).

\section{Algebraic description - The ideal-divisor-correspondence}

\noindent 1. {\it The affine case.} We remark that there is a
one-to-one correspondence between ideals of the quotient ring $\
{\CC} [x,y] /I_C$ of polynomial functions on the affine curve $C$
and the ideals of $\ {\CC} [x,y]$ with $I\supset I_C$ where $I_C
=(F(x,y,1))$. $\ {\CC} [x,y] /I_C$ is one-dimensional and a
Dedekind ring, cf. \cite{At}.

Now let $I$ be an ideal with $I\supset I_C$. Then $I$ is
zerodimensional. Because $\ {\CC} [x,y] /I_C$ is a Dedekind ring
we have the unique primary decomposition
\[ I = \bigcap (I_{P_i}^{m_i}, I_C)  \]
where $P_i=(x_i,y_i)$ and $I_{P_i}=((x-x_i),(y-y_i),I_C)$, cf.
\cite{At}. We associate to $I$ the effective divisor
\[D_I:=\sum m_i P_i\  .  \]
Conversely, let $D:=\sum m_i P_i$ be an effective divisor of
finite points. Then we associate to $D$ the ideal
\[ I_D:=\bigcap (I_{P_i}^{m_i}, I_C)\ .  \]
Therefore we have a one-to-one correspondence of effective
divisors and ideals with $I\supset I_C$.

Furthermore we have the following Lemma.
\begin{Le}
Let $D,D'$ be effective divisors of finite points and let $f\in
I_D$. Then we have $(f)\geq D$ and
\[
(I_D I_{D'},I_C) = I_{D+D'}\ \ \ \ \ \ \ \ {\rm and} \ \ \ \ \ \ \
\ \ \ \ \ ((f),I_C):I_{D}  = I_{(f)-D}\ .
\]
%(ii) Let $I,I'$ be two ideals with $I,I'\supset I_C$ and $f\in I$.
%Then we have $(f)\geq D_I$ and
%\[
%D_{(I \cdot\, I',I_C) }=D_I+D_{I'}\ \ \ \ \ \ \ \ {\rm and} \ \ \
%\ \ \ \ \ \ \ \ \ D_{(((f),I_C):I,I_C)}=(f)-D_{I}\ .
%\]
Here $(f)$ denotes both, the ideal generated by $f(x,y)$ and the
divisor of the analytic continuation of the meromorphic function
$f(x,y)$.
\end{Le}

\noindent {\it Proof.} The proof follows from
the above primary decomposition and the relations \\
$I_{P}^{p} \cap I_{Q}^{q} =
 I_{P}^{p}I_{Q}^{q}\ {\rm mod}\ I_C$
for $P\neq Q$ and $I_{P}^{p} \cap I_{P}^{q}=
 I_{P}^{\max\ (p,q)}\ {\rm mod}\ I_C$.
$\bullet$

\noindent 2. {\it The projective case.} We suppose that
$[0:1:0]\notin C$. Now let $D=D_e+D_\infty $ be an effective
divisor with the finite part $D_e$ and the infinite part $D_\infty
$. We define
\[ {\bf I}^h_D:=H_z(I_{D_e})\cap H_x(I_{D_\infty }) \]
where $I_{D_e}\subset\ {\CC} [x,y]$,
      $I_{D_\infty}\subset\ {\CC} [y,z]$
and $H_z$, $H_x$ are the homogenisations with respect to $z,x$,
respectively.

Conversely, let $\bf I$ be a homogeneous ideal of ${\CC} [x,y,z]$
containing the curve ideal, i.e. ${\bf I}\supset (F(x,y,z))$. Then
we define
\[ D_{\bf I}:=D_{A_z({\bf I})} +D_{(A_x({\bf I}),\ z)} \]
where $A_x, A_z$ are the affinisations with respect to $z,x$.

We have $D_{{\bf I}_D}=D$. We remark that for ideals ${\bf
I}_1,{\bf I}_2$ with ${\bf I}_1,{\bf I}_2\supset (F(x,y,z))$ the
relations ${\bf I}_{(D_{\bf I})}={\bf I}$ and ${\bf I}_{D_1} {\bf
I}_{D_2} = {\bf I}_{D_1+D_2} $ are in general not valid.

Now we define a product $\odot$ for homogeneous ideals ${\bf
I}_1,{\bf I}_2$ which corresponds to the addition of divisors. We
form the ideal product of the corresponding affine ideals. In
order to include infinite points we consider products with respect
to two affinisations with $z,x=1$. The intersection of the
corresponding homogenisations will contain all curves with
intersection divisor $\geq D_{{\bf I}_1}+D_{{\bf I}_2}$. I.e. we
define
\[
{\bf I}_1 \odot {\bf I}_2 := H_z((A_z({\bf I}_1) \cdot A_z({\bf
I}_2),I_C))\bigcap H_x((A_x({\bf I}_1) \cdot A_x({\bf
I}_2),I_C^x))
\]
where $I^x := A_x ( H_z(I))$ for ideals $I\subset {\CC} [x,y]$.

A generalized ideal quotient $\oslash$ is defined by
\[
(G)\oslash {\bf I} := H_z((A_z((G)),I_C  ) : A_z({\bf I}))\bigcap
                      H_x((A_x((G)),I_C^x) : A_x({\bf I})).
\]

\begin{Sa}
\noindent Let $D,D'$ be effective divisors and let $G\in {\bf
I}_D$. Then we have $(G)\geq D$ and
\[
{\bf I}^h_{D+D'}= {\bf I}^h_{D} \odot {\bf I}^h_{D'} \ \ \ \ \ \ \
{\rm and}\ \ \ \ \ \ \ {\bf I}^h_{(G)-D}= {\bf I}^h_{(G)} \oslash
{\bf I}^h_{D} .  \]
%\noindent (ii) Let ${\bf I,I'}$ be two ideals with ${\bf
%I,I'}\supset (F(x,y,z))$ and $G\in {\bf I}$. Then we have $(G)\geq
%D_{\rm I}$ and
%\[
%D_{{\bf I}\odot {\bf I}'}=D_{\bf I}+D_{{\bf I}'}\ \ \ \ \ \ \ \
%{\rm and}\ \ \ \ \ \ \ D_{(G) \oslash {\bf I}}= (G)-D_{{\bf I}}\ .
%\]
\end{Sa}

\noindent {\it Proof.} The proof follows from the fact, that the
left and the right side are equal to the homogeneous ideal of all
curves with intersection divisor $\geq D_{{\bf I}_1}+D_{{\bf
I}_2}$ and $\geq (G)-D_{{\bf I}}$, respectively. $\bullet$

\section{Reduction and the group law}

\noindent Let $D=D^+ - D^-$ with $\deg  D^+ =\deg D^- =s $ be a
divisor of degree zero and let ${\bf I^+}={\bf I}_{D^+}^h$, ${\bf
I^-}={\bf I}_{D^-}^h$. Furthermore we consider the homogeneous
ideals
\[
{\bf I}_r:={\bf I}^h_{rP_0}
\]
of forms with an $r$-fold point $P_0$. We choose $m$ such that
$b_m\geq s+g$. Then we choose an arbitrary element $G$ in ${\bf
I^+} \odot {\bf I}_{b_m-s} = {\bf I}^h_{D^+ +(b_m-s)P_0}$ of
degree $m$ and we form ${\bf J}=(G)\oslash{\bf ({\bf I}^+} \odot
{\bf I}_{ b_m-s })$. Then we determine the number $\alpha$ such
that $({\bf J}\odot {\bf I}^-)\cap {\bf I}_{b_m-s+g+\alpha }$
contains exactly one Groebner basis element $G'$ of degree $m$
with respect to a degree order. We form ${\bf I}_{red}=(G')\oslash
(({\bf J} \odot {\bf I}^- )\cap  {\bf I}_{b_m-s+g+\alpha } )= {\bf
I}_{red}={\bf I}^h_{S_1+\cdots S_{g-\alpha}}$ where $S_1+\cdots
S_{g-\alpha} -(g-\alpha)P_0$ is the reduced divisor of $D$.

\noindent {\it Remark}: Given ${\bf I}^+$, ${\bf I}^-$, $t$ one
can carry out the determination of ${\bf I^+} \odot {\bf
I}_{b_m-s}$, ${\bf J}=(G)\oslash{\bf ({\bf I}^+}\odot {\bf I}_{
b_m-s })$, $({\bf J}\odot {\bf I}^-)\cap {\bf I}_r$ $(r\geq
b_m-s+g )$ and $(G')\oslash (({\bf J}\odot {\bf I}^-)\cap {\bf
I}_{b_m -s+g+ \alpha })$ by Groebner basis calculations, cf.
\cite{Be}.

Now we can describe the group law.
In order to add the divisors
$D_1=P_1+\ \cdots\ P_g-gP_0$ and $D_2=P_{g+1}+\ \cdots\
+P_{2g}-gP_0$ we apply the above construction to
$D^+=P_1+\ \cdots\ +P_{2g}$ and $D^-=2gP_0$.
In order to determine the ideal for $-D$ of $D=P_1+\ \cdots\
+P_g-gP_0$ (where $P_1,\ \cdots\ ,P_g = P_0$ is allowed)
we apply the above construction to $D^+=gP_0$ and $D^-=P_1+\
\cdots\ +P_g$.

\section{An Example}

\noindent We consider the 4-curve $C$ with $x^4+y^4=2z^4$ with $g=3$.
Let $P_0=(1,1,1)$ and $P_1=(1,-1,1)$. We reduce the divisor
\[
D=6P_1-6P_0.
\]
We choose $m=3$. It is sufficient to consider affine ideals. ($C$
has no infinite intersections with all curves occuring during the
calculation.) We obtain by Groebner basis calculations with
respect to a lexicographic order
\[
I^+=\{1 + 6y + 15y^2 + 20y^3 + 15y^4 + 6y^5 + y^6,
     102 - x + 524y + 1092y^2 + 1141y^3 +
     598y^4 + 126y^5 \},
\]
\[
I_3=\{-1 + 3y - 3y^2 + y^3, -1 - x + 5y - 3y^2\},
\]
\[
(I^+I_3 ,I_C) =\{-1 - 3y + 8y^3 + 6y^4 - 6y^5 - 8y^6 + 3y^8 + y^9,\]
\[
 1289 - 128x + 2492y - 2016y^2 - 7476y^3 - 806y^4 + 7476y^5 + 3080y^6 -
2492y^7 -
  1419y^8\}.
\]
With respect to a degree order we find an element of degree 3 in
$(I^+I_3,I_C)$
\[
f=2612 - 3078x + 378x^2 + 281x^3 + 478y - 1912xy + 1195x^2y - 1286y^2
+ 1093xy^2 +
 239y^3.
\]
Furthermore
\[
J:=((f),I_C):(I^+I_3) =\{698405268857 + 635735348837y -
10585774871y^2 + 108619669441y^3,\]
\[
 268707776349 + 254869376165x + 103445986821y + 108619669441y^2 \},
\]
\[
- 66326x^2y -
 735502y^2 + 70382xy^2 + 325163y^3,
\]
\[
I_6=\{1 - 6y + 15y^2 - 20y^3 + 15y^4 - 6y^5 + y^6, -102 + x + 524y -
1092y^2 + 1141y^3 -
  598y^4 + 126y^5\},
\]
{\tiny
\[
(JI_6 ,I_C) =\{698405268857 - 3554696264305y + 6651081164962y^2 -
4259940825918y^3 -
  3049132583596y^4 + 7186609158448y^5 - 5447186836050y^6 \]
  \[
  + 2328545039678y^7 - 662303791517y^8 + 108619669441y^9,\]
\[
 -93764434230515570410626334329726164533810553558278390921760183 +
  2185848162978035880543485328036429801750462217550638009464820x \]\[+
  390604925102813740313936712629510601523986718107312688345228286y -
  553599197596387095117790236962526646313212938826438985782932370y^2\]\[ +
  109584143933719697324470424983767313360642548080131427940056882y^3 +
  488838345917420338620885461620364216394401215438274344792992080y^4 \]\[-
  548212607851251833242108578943280138884947131311429528021917262y^5 +
  272507862737893198250637137223325900667162331316521096907832914y^6\]\[ -
  85848455024362109011538768963037701794199996136791647052272470y^7 +
  17703568847691597391590697413566189778227344673148355783307303y^8\}.\]}
There is an element of degree 3 in $(JI_6,I_C)$
\[
g=683086 - 414993x - 636078x^2 + 356233x^3 - 259643y + 677678xy. \]
Finally we obtain
\[
I_{red}= ((g),I_C):(JI_6 ) = \{94544281343 + 377260313207y +
408415639297y^2 + 134215744153y^3,\]\[
 -53515118937 + 13173978910x - 225487128300y - 134215744153y^2\}.
\]
Because the minimal element of $I_{red}$
\[
-53515118937 + 13173978910x - 225487128300y - 134215744153y^2
\]
with respect to a degree order is of degree $2>n-3$ the ideal
$I_{red}$ is reduced. The ideal $I_{red}$ corresponds to the
reduced divisor
\[
 (-0.82409 - 0.62806i,
  -1.31975 + 0.06425i,1)\]\[+
 (-0.82409 + 0.62806i,
  -1.31975 - 0.06425i,1)+
 (-1.18524,- 0.40347,1)-3(1,1,1)\sim D.
\]

\section{Curves with simple double points }

The above construction applies analogously to curves with simple
singularities. Here we consider the case of $n$-curves $C$ with
$d$ finite simple double points $ {\cal D}_1,\ \cdots\ ,{\cal D}_d \neq P_0 $
with $F_x=F_y=0,\ F_{xx} F_{yy}-F_{xy}^2\neq 0$. Furthermore we
suppose $[0,1,0]\notin C$. To every point ${\cal D}_i$ correspond two
points ${\cal D}_i^+$, ${\cal D}_i^-$ on the Riemann surface of $C$ of genus
$g=\frac{(n-1)(n-2)}{2}-d$. Let
\[
\Delta:={\cal D}_1^+ +{\cal D}_1^- +\ \cdots\ +{\cal D}_d^++{\cal D}_d^-
\]
be the double point divisor of $C$. Now consider the divisor
\[
D=D^+ - D^- =P_1+\ \cdots\ +P_s-Q_1-\ \cdots\  -Q_s
\]
of degree zero. We consider an $m$-curve with a polynomial
$G(x,y,z)\neq 0\ {\rm mod}\ F$ with
\[
s+g+d\leq b_m:=\left\{  \begin{array}{ccc}
\frac{m(m+3)}{2} & {\rm for} & m<n \\
mn-\frac{(n-1)(n-2)}{2} & {\rm for} & m\geq n
\end{array} \right.\
\]
through the $s$ points of $D^+$, ${\cal D}_1,\ \cdots\ ,{\cal D}_d$ and
$(b_m-s-d)P_0$. We have \[ mn-s-2d-(b_m-s-d)\ =\ mn-b_m-d\ =\
\frac{(n-1)(n-2)}{2}-d\ =\ g \] remaining intersections $R_1,\
\cdots\ ,R_g$, i.e.
\[
D^+ +\Delta +(b_m -s-d )P_0+R_1+\ \cdots\ +R_g \sim m D_\infty .
\]
Then we consider an $m$-curve $G'(x,y,z)\neq 0 \rm\ mod\ \it F$
through the $s$ points of $D^-$ and through $R_1,\ \cdots\ ,R_g,
{\cal D}_1,\ \cdots\ ,{\cal D}_d,(b_m-s-d-g)P_0$. We require a maximal
additional contact $\alpha$ at $P_0$. Let $S_1,\ \cdots\
,S_{g-\alpha}$ be the remaining intersections not equal to $P_0$.
We have
{\small  \[
D^- + \Delta + R_1+ \cdots +R_g +( b_m-s-d-g+\alpha )P_0 +S_1+
\cdots +S_{g-\alpha}
 \sim m D_\infty .
\]}
It follows
\[
D^+ - D^-  \sim \overline{D}:= S_1+\ \cdots\ +S_{g-\alpha} -
(g-\alpha)P_0 .
\]
Analogously to Proposition 4 one shows that $\overline{D}$ is the
reduced divisor of $D$.

For an algebraic description of the reduction let ${\bf I}_\Delta
$ be the ideal of all adjoint curves through the ${\cal D}_1,\ \cdots\
,{\cal D}_d$. For a finite simple double point ${\cal D}_i$ we have two taylor
series expansions
\[
y-y_i= a^\pm_1(x-x_i)+ a^\pm_2(x-x_i)^2+\cdots\
\]
for the two branches of ${\cal D}_i^\pm$ (or $x-x_i= b^\pm_1(y-y_i)+
b^\pm_2(y-y_i)^2+\cdots\  $, if $b_1^\pm=0$).

Then
{\small \[
{\bf I}_{m {\cal D}_i^\pm}:=H_z( ( (x-x_i,y-y_i)^m,
y-y_i-a^\pm_1(x-x_i)-
       \cdots -a^\pm_{m-1}(x-x_i)^{m-1} )  )
\]}
are the ideals of polynomials with an $m$-fold common point with
the branch of ${\cal D}_i^\pm $. Now let
\[  D=m_1P_1+\ \cdots\ +m_pP_p+\sum_i (m^+_i {\cal D}_i^+ + m^-_i {\cal D}_i^-)
\]
be a divisor with different ordinary points $P_1,\ \cdots\ ,P_p$.
We consider the ideal of curves with intersection divisor $ \geq
\Delta +D $
\[
{\bf I}_{\Delta +D}^h:= ({\bf I}_{m_1P_1}^{h}\cap \cdots \cap {\bf
I}_{m_pP_p}^{h}) \cap \bigcap_i ( {\bf I}_{(m_i^+ +1) {\cal D}_i^+}
            \cap {\bf I}_{(m_i^- +1) {\cal D}_i^-} ).
\]
%where ${\bf I}_{P}^{k}:=H_z(I_P^k,I_C)
%\cap H_x((I_P^x)^k,I_C^x) $.
Now we define an ideal product $\odot_\Delta $ for ideals $ {\bf
I}_1={\bf I}_{\Delta +D_1}^h, {\bf I}_2={\bf I}_{\Delta +D_2}^h $
which corresponds to the addition of the divisors $D_1$ and $D_2$.

{\footnotesize
\[
{\bf I}_1 \odot_\Delta {\bf I}_2 := H_z(  ( A_z({\bf I}_1) \cdot
A_z({\bf I}_2),I_C ) : I_\Delta     )
    \bigcap
H_x(  ( A_x({\bf I}_1) \cdot A_x({\bf I}_2),I_C^x  ) : A_x({\bf
I}_\Delta ) )
% \bigcap
 \]}
%{\footnotesize
%\[
%H_z(\ ( \cap_{i,\varepsilon
%= \pm }
%     (A_z({\bf I}_1) \cdot A_z({\bf I}_2),J_{ D_i^{ \varepsilon } ,\ t } ) : I_\Delta ,
%       J_{ D_i^{ \varepsilon } ,\ t }
%    ) \ )
%    \bigcap
% H_x(\ ( \cap_{i,\varepsilon = \pm }
%     (A_x({\bf I}_1) \cdot A_x({\bf I}_2),J^x_{ D_i^{ \varepsilon },\ t } ) : I_\Delta^x ,
%       J^x_{ D_i^{ \varepsilon },\ t }
%    )\  )
%\]}
\noindent where $ I_C^x := A_x ( {F} ) $.
%and where $t$ is sufficient large.
${\bf I}_1 \odot_\Delta {\bf I}_2$ contains all curves with
intersection divisor $\geq \Delta +D_1+D_2$.

A generalized ideal quotient is given by {\footnotesize
\[
(G)\oslash_\Delta {\bf I} := H_z(  (A_z( (G) )\cdot I_\Delta ,I_C)
: A_z({\bf I}) ) \bigcap H_x(  (A_x( (G) )\cdot A_x({\bf I}_\Delta
),I_C^x) : A_x({\bf I}))
%\bigcap
\]}
%\[
%H_z(\  \cap_{i,\varepsilon =\pm }
%    (A_z( (G) )\cdot I_\Delta,J_{ D_i^{ \varepsilon } ,\ t }) :
%    A_z({\bf I}) \  )
%    \bigcap
%H_x(\  \cap_{i,\varepsilon =\pm }
%    (A_x( (G) )\cdot I_\Delta^x,J^x_{ D_i^{ \varepsilon },\ t }) :
%    A_x({\bf I}) \  )
%    \   .
%\]}
%\[
%{\rm I}_1 \odot_\Delta {\rm I}_2 :=
%\cap_{i=1}^d
%\left(
%({\rm I}_1 \odot {\rm I}_2,{\bf I}_{m D_i^+} )
%\cap
%({\rm I}_1 \odot {\rm I}_2,{\bf I}_{m D_i^-} )
%\right)
%\]
%where the $m_{\varepsilon_i}$ are sufficient large.
\begin{Sa}
Let $D,D'$ be effective divisors and let
$G\in {\bf I}_{\Delta +D}$. Then we have $(G)\geq \Delta +D$ and
\[
{\bf I}^h_{\Delta +D+D'}= {\bf I}^h_{\Delta +D} \odot_\Delta {\bf
I}^h_{\Delta +D'} \ \ \ \ \ \ \ {\rm and}\ \ \ \ \ \ \ {\bf
I}^h_{(G)-D}= {\bf I}^h_{(G)} \oslash_\Delta {\bf
I}^h_{\Delta +D} .  \]
\end{Sa}

\noindent {\it Proof.}
 The proof follows from the fact,
that the left and the right side are equal to the homogeneous
ideal of all curves with intersection divisor $\geq \Delta +
D_{{\bf I}_1}+D_{{\bf I}_2}$ and $\geq (G)-D_{{\bf I}}$,
respectively. $\bullet$

Now let $D=D^+ - D^-$ with $D^+,D^-\geq 0$, $\deg  D^+ =\deg D^-
=s $ be a divisor of degree zero. Furthermore let ${\bf I^+}:={\bf
I}_{\Delta +D^+ }^h  $, ${\bf I^-}:={\bf I}_{\Delta +D^- }^h  $,
${\bf I}_r:={\bf I}^h_{rP_0}$. We choose $m$ such that $b_m\geq
s+d+g$. Then we choose an arbitrary element $G$ in ${\bf I^+}
\odot {\bf I}_{b_m-s-d} = {\bf I}^h_{\Delta  + D^+
+(b_m-s-d)P_0}$ of degree $m$ and we form ${\bf J}=(G)
\oslash_\Delta ({\bf I^+} \odot {\bf I}_{b_m-s-d}) ={\bf
I}^h_{\Delta +R_1+\ \cdots\ +R_g}$. We determine the number $r\geq
b_m-s-d-g$ such that $({\bf J}\odot_\Delta {\bf I}^-)\cap {\bf
I}_r$ contains exactly one Groebner basis element $G'$ of degree
$m$ with respect to a degree order. We form ${\bf
I}_{red}=(G')\oslash_\Delta ( ({\bf J} \odot_\Delta {\bf I}^- )
\cap {\bf I}_r )$. We obtain ${\bf I}_{red}={\bf I}^h_{\Delta
+P_1+\cdots +P_t}$ with $t=g-r+(b_m-s-d-g)$ where $P_1+\cdots\
+P_t-tP_0$ is the reduced divisor for $D$.

\section{An Example}

We consider the hyperelleptic $4$-curve $C$ with $x^{4}-y^{4}=30xyz^2$
with $g=2$.
$C$ has one simple double point ${\cal D}_0^\pm =    (0,0,1)$ with the
Taylor series expansions
\[
y=   \frac{x^3}{30} - \frac{x^{11}}{24300000} +\ \cdots ,\ \ \ \ \ \ \
x= - \frac{y^3}{30} + \frac{y^{11}}{24300000} -\ \cdots .
\]
We choose $P_0:=(1,1,0)$ and we consider the points
$P_1=(4,2,1)$, $P_2=(1,-1,0)$.
We apply the reduction to ideals ${\bf I}^+,{\bf I}^-$, which correspond to
the divisor $2{\cal D}_0^+-P_1-P_2$.

We have
\[  {\bf I}^+ = {\bf I}^h_{\Delta+2{\cal D}_0^+}  = (x^3,y) ,   \]
and
\[  {\bf I}^- = {\bf I}^h_{\Delta+P_1+P_2} =
(x z - 2 y z, -x y - y^2 + 6 y z, -x^2 + y^2 + 6 y z,
    y^2 z - 2 y z^2).    \]
We form
\[  {\bf I}^h_{\Delta+2{\cal D}_0^+ +2P_0} = {\bf I}^+ \odot (x-y,z^2)  =
( -yx + y^2, y z^2, x^3 -x^2 y, x^3 z^2). \]
We choose $G=-xy+y^2$. We obtain the quotient
\[ {\bf J}=(G)
\oslash_\Delta {\bf I}^h_{\Delta+2{\cal D}_0^+ +2P_0} =
(x^2,xy,y^2) =
{\bf I}^h_{\Delta +R_1 +R_2}.    \]
We remark that $R_1+R_2={\cal D}_0^++{\cal D}_0^-$. Now we form
\[
{\bf I}^h_{\Delta+R_1+R_2+P_1+P_2} = {\bf I}^- \odot_\Delta J=
 (-x^2 + x y + 2 y^2, x^2z - 2 xyz,
  x^3 +x^2 y - 6x^2 z, x^3z - 4 x^2z^2). \]
We choose $G'=x^2-xy-2y^2$. Then we obtain the quotient
\[
(G') \oslash_\Delta {\bf I}^h_{\Delta+R_1+R_2+P_1+P_2} =
(x z - 2 y z, -x y - y^2 - 6 y z, x^2 - y^2 + 6 y z,
    y^2 z + 2 y z^2) =
{\bf I}^h_{\Delta +S_1 +S_2}={\bf I}_{red}.
\]
We obtained ${\bf I}_{red}$ from ${\bf I}^+,{\bf I}^-$ by
rational operations. We remark that
$S_1+S_2=(-4,-2,1)+P_2$.

\section{Hyperelliptic curves}

\noindent
In this section we present an algorithm for hyperelliptic curves $C$
of genus $g$ in the standard form
\[   y^2=a(x-x_1)(x-x_2) \cdots (x-x_{2g+1})=:h(x)  \]
with different $x_i$, $a\neq 0$.
The projectivisation has a single singular (nonsimple)
infinite point $P_\infty =(0,1,0)$.
We choose $P_0:=P_\infty$.
Hyperelliptic curves have an involution
$x\rightarrow x$, $y\rightarrow -y$.
Let $\overline{I}$ and $\overline{D}$ be the image
of the ideal $I$ and the divisor $D$
with respect to this involution. Then we have
$\overline{(x_1,y_1)+\cdots +(x_n,y_n)}=
(x_1,-y_1)+\cdots +(x_n,-y_n)$,\\
$\overline{(x,y)}+(x,y)\sim 2 P_0$
and $I_{\overline{D}}=\overline{I_D}$.

Because $C$ has only one infinite point
it is sufficient to consider affine ideals.
However, the selection of the interpolating curves
requires a modification.

Let $D=D^+ - D^- \sim D^+ +\overline{D^-}-2\deg(D^-)P_0$
be a divisor of finite points.
We have $J:=I_{D^+ +\overline{D^-}}=I_{D^+}\overline{I_{D^-}}$.
We replace the degree order by the weighted degreelexicographic
order $\deg_{2g+1,2}$
with $\deg_{2g+1,2}(x^ay^b):=(2g+1)a+2b$ and $y>x$.
Now we determine the minimal element $f=p(x)+q(x)y$ of $J$
with respect to this order.
Then $C$ and the curve $f(x,y)=0$ have $\deg_{2g+1,2} (f)$
finite intersections
with $p^2(x)-h(x)q^2(x)=0$ and $y=p(x)/q(x)$.
Let $(x_i,y_i)$, $i=1,\cdots ,t:=\deg_{2g+1,2}(f)-\deg(D^+)-\deg(D^-)$
be the remaining finite intersections.
It follows
\[  (f)= D^+ +\overline{ D^-} + (x_1,y_1)+\cdots +(x_{t},y_{t})
 - \deg_{2g+1,2}(f) P_0    ,  \]
i.e.
\[   D-(\deg(D^+)-\deg(D^-))P_0
\sim (x_1,-y_1)+\cdots +(x_{t}, - y_{t}) -tP_0  =:D_1-tP_0. \]
\begin{Le}
$D_1-tP_0$ is the reduced divisor for $D-\deg(D)P_0$.
\end{Le}

\noindent
{\it Proof:} The proof follows from the
fact that $t:=\deg_{2g+1,2}(f)-\deg(D^+)-\deg(D^-)$
is minimal if $\deg_{2g+1,2}(f)$ is minimal.
$\bullet$

The divisor
$D_1$ corresponds to the  ideal
\[  I_{red} = I_{D_1}=\overline{(f,I_C):(I_{D^+}\overline{I_{D^-}})}.  \]

\noindent
{\it Remark:}
Contrarily to Cantors algorithm we can describe our algorithm
by this single formula. Our algorithm uses a reduction function
of the general form
$f=p(x)+yq(x)$. In contrast Cantors algorithm uses reduction functions
of the special form $y-p(x)$ several times
(cf. the example below).

\section{Picard curves}

\noindent
In this section we present an algorithm for Picard curves $C$
\[   y^3=a(x-x_1)(x-x_2)(x-x_3)(x-x_4)=:h(x)  \]
with four different $x_i$, $a\neq 0$.
The projectivisation has
a single 4-fold infinite point $P_\infty =(0,1,0)$.
We choose $P_0:=P_\infty$.
$C$ has only one infinite point and
it is sufficient to consider affine ideals.
Similarly to the case of hyperelliptic curves
we replace the degree order by the weighted degreelexicographic order
with $deg_{4,3}:=(x^ay^b):=4a+3b$ and $y>x$.

Let $D=D^+ - D^- $
be a divisor of finite points.
We determine the minimal element
$f=p(x)+q(x)y+r(x)y^2$ of $I_{D^+}$
with respect to the above order. We have
\[ (f) = D^+ + R_1 + \cdots + R_t  - \deg_{4,3} (f) P_0 \]
with $t:=\deg_{4,3}(f)-\deg(D^+)$
remaining finite points $D'=R_1 + \cdots + R_t$
whose $x$-coordinates are zeros of the polynomial
\[ \left| \begin{array}{ccc}
p & q & r \\
hr & p & q \\
hq & hr & p
\end{array} \right| =p^3+q^3h+r^3h^2-3pqrh  \]
of degree $t$.
We have
\[   I_{D'} = (f,I_C):I_{D^+}. \]
We determine the minimal element $g$ of $I_{D'+D^-}=I_{D'}I_{D^-}$
with respect to the above order. We have
\[ (g) = D'+D^- + S_1 + \cdots + S_q  - \deg_{4,3} (g)P_0 \]
with $q:=\deg_{4,3}(g)-\deg(D')-\deg(D^-)$
remaining finite points $D''=S_1 + \cdots + S_q$.
It follows
\[
D''-qP_0 \sim -D'-D^- +(\deg_{4,3}(f)-q)P_0
\sim D-(\deg(D^+)-\deg(D^-))P_0.
\]

\begin{Le}
$D''-q P_0$ is the reduced divisor for $D- \deg(D) P_0$.
\end{Le}

\noindent
{\it Proof:} The proof follows from the
fact that $q$ is minimal if $\deg_{4,3}(g)$ is minimal.
$\bullet$

We have
\[  I_{red} = I_{D''} = (g,I_C):I_{D'+D^-}. \]
An algorithm for the Jacobian group of this curve is discussed
in \cite{Es,FO}.
Contrarily to our algorithm these algorithms
use the concrete structure of the curve and
require the distinction of many different cases.
Our algorithm for Picard curves
has a straightforward generalization to superelliptic
curves
\[ y^m=a(x-x_1)\cdots (x-x_n)  \]
with different $x_i$, $a\neq 0$ and ${\rm lcd} (m,n)=1$.

\section{The case of characteristic $p$}

\noindent
In this last section we make some remarks about the case of
characteristic $p$.
Using the theory of \cite{Mo} one can show that our algorithm
has an analogue if we replace\ $\CC$ by a field
$k$ of characteristic $p$.
$k$ has the algebraic closure
$\overline{k}$. The divisor group $Div(C)$ is the free Abelian group
consisting of formal finite sums $\sum_{P\in C(\overline{k})}m_P P$
with $m_P\in \ZZ$.

A divisor is defined over $k$ if it is fixed by the natural Galois
action of ${\rm Gal} (\overline{k},k)$.
The divisors defined over $k$ form the subgroup $Div_k(C)$.
Analogously one defines the group $Div^0_k(C)$.
Principal divisors are defined as zeros and poles of rational
functions $\frac{G(x,y,z)}{H(x,y,z)}$
where $G,H$ $(H\notin (F(x,y,z))$  are homogeneous polynomials
of equal degree
with coefficients in $k$.
They form a subgroup $Div^P_k(C)\subset Div^0_k(C)$.
We define $Jac_k(C):=Div^0_k(C) / Div^P_k(C)$.
Furthermore we can define the analog notion of an reduced divisor.

Let $D$ be an element of $Div_k(C)$ without infinite points.
Then the polynomials $p(x,y)\in k[x,y]$
with $(p)\geq D$ form an ideal $I_D\subset k[x,y]$.
One can show that there is a
one-to-one correspondence between ideals $I$ of $k[x,y]$ with
$I\supset I_C$ and divisors $D$ of $Div_k(C)$ of finite points.

Now let $C$ be a hyperelliptic curve or a Picard curve.
The affine part of the above hyperelliptic curves is smooth
for $p\neq 2$ and the above Picard curves are smooth for
$p\neq 3$.
This case is interesting in view of applications in cryptography.
Let
$I^+,I^-$ be two ideals of $k[x,y]$ with divisors $D^+,D^-$.
We apply the corresponding algorithm of the two previous sections
to $I^+,I^-$. Because all operations are rational
we obtain an ideal $I_{red}$ of $k[x,y]$.
We have an analogue of Lemma 9,10. Therefore
$I_{red}$ corresponds to a reduced divisor.

\noindent
{\it Example:} Let $C$ be the hyperelliptic curve
$y^2=(x-3)(x-2)(x-1)x(x+1)(x+2)(x+3)=:h(x)$ for $k=F_{17}$
and let
$D_1=(4,5)+(5,8) +(6,4)$,
$D_2=(7,5)+(10,3)+(11,1)$.
We have
$I_{D_1}=(x^3+2x^2+6x+16,5x^2+9x+8-y)=:(a_1,b_1-y)$
and
$I_{D_2}=(x^3+6x^2+2x+12,11x^2+5x+9-y)=:(a_2,b_2-y)$.

\noindent
1. \it Cantors algorithm: \rm
We have $s_1a_1+s_2a_2=1$ with certain
$s_1$,
$s_2$.
We obtain the composed ideal $I_{D_1}I_{D_2}$ with the basis
$(a,b-y):=(a_1a_2,(s_1a_1b_2+s_2a_2b_1-y) {\rm mod} a_1a_2)=
(x^6+8 x^5+3x^4+13x^2+2x+5,x^5+7x^4+2x^3+6x^2+5x+5-y )$.
The reduction process gives the ideals of equivalent divisors \\
$(a',b'-y):=(\frac{b^2-h}{a},(-b-y){\rm mod} a')=
 (x^4+6x^3+2x^2+5x+5, 6x^3+x^2+5x-y )$ and \\
$(a'',b''-y):=(\frac{{b'}^2-h}{a'},(-b'-y){\rm mod} a'')=
(x^3+9x^2+3x,2x^2+13x-y)$.\\
The last ideal corresponds to the reduced divisor.

\noindent
2. \it Our algorithm: \rm
We obtain for $I_{D_1}I_{D_2}$ the Groebner basis
$(x^6+8 x^5+3x^4+13x^2+2x+5,11x^4+9x^3+2x^2+9x+y(x+1))=:(.,f)$
with respect to the weighted lexicographic order
and
$(f,y^2-h(x)):(I_{D_1}I_{D_2})= (x^3+9x^2+3x,2x^2+13x-y)$.


\begin{thebibliography}{9999999}

\bibitem[1]{Ar} S. Arita, {Algorithms for computations
                 in Jacobian group of $C_{ab}$ curve and their
                 application to discrete-log-based public key
                 cryptosystems}, in A. Odlyzko et al (eds.),
                 The mathematics of public key cryptography
                 Fields Institute, Toronto, 1999.

\bibitem[2]{At} M.F. Atiyah, I.G. Macdonald,
                 { Introduction to commutative algebra}
                 (Addison-Wesley, Reading, 1969).

\bibitem[3]{Be} Th. Becker, V. Weispfenning,
                 { Gr\"{o}bner bases} (Springer, New York, 1993).

\bibitem[4]{Ca} D.G. Cantor, { Computing in the Jacobian
                 of a hyperelliptic curve},
                 Mathematics of Computation
                 { 48} 177(1987) 95-101.

\bibitem[5]{Cl} A. Clebsch, P. Gordan,
                 Theorie der Abelschen Functionen,
                 Teubner, Leipzig : Teubner, 1866.

\bibitem[6]{Ei} D. Eisenbud, M. Green, J. Harris,
                 { Cayley-Bacharach theorems and conjectures},
                 Bull. Am. Math. Soc., New Ser. { 33}
                 (1996) 295-324.

\bibitem[7]{Es} J. Estrada Sarlabous, E. Reinaldo Barreiro,
                 J.A. Pineiro Barcelo, { On the Jacobian Varieties of
                 Picard Curves}, Math. Nachr. { 208}
                 (1999) 149-166.

\bibitem[8]{FO} S. Flon, R. Oyono,
                 Fast arithmetic on Jacobians of Picard curves,
                 LNCS { 2947} (2004) 55-68.


\bibitem[9]{Ga} S.D. Galbraith, S.M. Paulus, N.P. Smart,
                 { Arithmetic on superelliptic curves},
                 Math. Comput. { 71} (2002) 393-405.

\bibitem[10]{Har} R. Harasawa, J. Suzuki,
                 { Fast Jacobian group arithmetic on $C_{ab}$ curves},
                 Lect. Notes Comput. Sci., vol. 1838 (2002) 359-376.

\bibitem[11]{Ha} R. Hartshorne, { Algebraic geometry},
                 Springer, New York 1987.

\bibitem[12]{He} F. He{\ss}, { Computing Riemann-Roch spaces
                 in algebraic function fields and related topics},
                 J. Symb. Comput. { 33} (2002) 425-445.

\bibitem[13]{Hu} M.-D. Huang, D. Ierardi, { Efficient Algorithms
                 for the Riemann-Roch Problem and for Addition  in
                 the Jacobian of a Curve},
                 J. Symb. Comp. { 18} (1994) 519-539.

\bibitem[14]{K}  N. Koblitz,
                 { Algebraic aspects of cryptography.
                 With an appendix on hyperelliptic curves},
                 Springer, New York, 1999.

%\bibitem[Maz86]{Ma} B. Mazur, { Arithmetic on curves},
%                 Bull. AMS { 14}(1986), 207-259.

%\bibitem[13]{Ma} Khuri-Makdisi, K.,
%                 Linear algebra algorithms for divisors on an algebraic
%curve,
%                 preprint math.NT/0105182.

%\bibitem[14]{Mu1}Mumford, D., On the equations defining Abelian
%varieties,
%                 I-III,
%                 Invent. Math. 1(1966)287-354,
%                 ibid. 3(1967)75-135, 215-244.

\bibitem[15]{Mo} C. Moreno, Algebraic curves over finite fields,
                 Cambridge Univ. Pr., Cambridge, 1993.

\bibitem[16]{Mu} D. Mumford,
                 { Tata lectures on theta},
                 Birkh\"{a}user, Boston, 1994.

\bibitem[17]{Vo} E.J. Volcheck, { Computing in the Jacobian of a Plane
                 Algebraic Curve},
                 ANTS-I, Springer LNCS { 877} (1994) 221-233.
\end{thebibliography}
\end{document}